\documentclass[reqno]{amsart}
\usepackage{hyperref}

\begin{document}
\title[\hfilneg EJDE-2006/71\hfil Energy quantization]
{Energy quantization for Yamabe's problem in conformal dimension}
\author[F. Mahmoudi\hfil EJDE-2006/71\hfilneg]
{Fethi Mahmoudi}

\address{Fethi Mahmoudi \newline
Scuola Internazionale Superiore Di Studi Avanzati (Sissa)\\
Via Beirut 2-4, 34014 Trieste, Italy}
\email{mahmoudi@sissa.it}

\date{}
\subjclass[2000]{35B33, 46E30, 46L65}
\keywords{Critical exponents; Lorentz spaces; quantization phenomena}

\begin{abstract}
 Rivi\`ere \cite{T.R} proved an energy quantization for
 Yang-Mills fields defined on $n$-dimensional Riemannian manifolds,
 when $n$ is larger than the critical dimension 4. More precisely,
 he proved that the defect measure of a weakly converging sequence
 of Yang-Mills fields is quantized, provided the $W^{2,1}$ norm
 of their curvature is uniformly bounded. In the present paper,
 we prove a similar quantization phenomenon for the nonlinear
 elliptic equation
 \[
 - \Delta{u}= u  |u|^{4/(n-2)},
 \]
 in  a subset $\Omega$ of $\mathbb{R}^n$.
\end{abstract}

\maketitle \numberwithin{equation}{section}
\newtheorem{theorem}{Theorem}[section]
\newtheorem{lemma}[theorem]{Lemma}
\newtheorem{remark}[theorem]{Remark}

\section{Introduction}

Let $ \Omega  $ be an open subset of $\mathbb{R}^n$ with
$n\geq 3$. We consider the equation
\begin{equation}
 -\Delta{u}=u|{u}|^{4/(n-2)}   \qquad \hbox{ in}\quad   \Omega
\label{eq1.1}
 \end{equation}
We will say that $u$ is a weak solution of \eqref{eq1.1} in
$\Omega$, if, for all  $ \Phi \in C^{\infty}({\Omega})$  with
compact support in
 $ \Omega $, we have
\begin{equation}
-\int_{\Omega}{\Delta{\Phi}(x) u(x)} dx = \int_{\Omega}{\Phi}(x)
 u(x)|{u(x)}|^{4/(n-2)}dx  \label{eq1.2}
 \end{equation}
 If in addition $u$ satisfies
\begin{equation}
\int_{\Omega}\Big[{\frac{\partial{u}}{\partial{x_{i}}}}
{\frac{\partial{u}}{\partial{x_{j}}}}
{\frac{\partial{\Phi^{j}}}{\partial{x_{i}}}}-\frac{1}{2}|\nabla{u}|^{2}
{\frac{\partial{\Phi^{i}}}{\partial{x_{i}}}}
+\frac{n-2}{2n}|u|^{2n/(n-2)}
{\frac{\partial{\Phi^{i}}}{\partial{x_{i}}}}\Big]dx  = 0
\label{eq1.3}
\end{equation}
for any $\Phi=(\Phi^{1},\Phi^{2}\dots ,\Phi^n)$ $\in
C^{\infty}(\Omega)$ with compact support in $\Omega$, we say that
u is stationary.
In other words, a weak solution $u$ in
$\textbf{H}^{1}(\Omega)\cap \textbf{L}^{2n/(n-2)}(\Omega)$
of \eqref{eq1.1}  is stationary if the functional $E$ defined
by
\[
E(u)=\frac{1}{2}\int_{\Omega}|\nabla{u}|^{2}
+\frac{n-2}{2n}\int_{\Omega}|u|^{2n/(n-2)}
\]
is stationary with respect to  domain variations, i.e.
\[
\frac{d}{dt}(E(u_{t}))|_{t=0}=0
\]
 where $u_{t}(x)=u(x+t\Phi)$.  It is easy to verify that a smooth
solution is stationary.

 In this paper we prove a monotonicity formula for stationary weak
solution $u$ in
$\textbf{H}^{1}(\Omega)\cap \textbf{L}^{2n/(n-2)}(\Omega)$
of \eqref{eq1.1} by a similar idea as in \cite{F.P}. More precisely
we have the following result.

\begin{lemma} Suppose that
$u \in \mathbf{L}^{2n/(n-2)}(\Omega ) \cap \mathbf H^1 (\Omega)$ is a
stationary weak solution of \eqref{eq1.1}. Consider the function
\[
E_u (x,r) =\int_{B(x,r)}|u|^{2n/(n-2)}dy +
\frac{d}{dr}\int_{\partial B(x,r)}u^2 ds +
r^{-1}\int_{B(x,r)}u^2 ds.
\]
Then  $r \mapsto E_u(x,r) $ is positive, nondecreasing and continuous.
 \label{lem1.1}
\end{lemma}

This monotonicity formula together with ideas which go back to the
work of  Schoen \cite{Sch}, allowed  to prove the following result.

\begin{theorem}
 There exists $ \varepsilon>0 $ and $r_{0}>0 $ depend only on $n$ such that,
for any smooth solution  $u\in \textbf{H}^{1}(\Omega)\cap
\textbf{L}^{2n/(n-2)}(\Omega)$ of \eqref{eq1.1}, we have:
For any $x_{0}\in\Omega$, if
\[\int_{B(x_{0},r_{0})}|\nabla{u}|^{2}+|u|^{2n/(n-2)} \leq
\varepsilon,
\]
then
\[\|u\|_{\textbf{L}^{\infty}(B_{\frac{r}{2}}(x_{0}))}\leq\frac{C(\varepsilon)}{r^{(n-2)/2}}
\quad\hbox{for any}  \quad r<r_{0},
\]
where $B_{\frac{r}{2}}(x_{0})$ is the ball centered at $ x_{0}$ with
radius $\frac{r}{2}$  , and $C(\varepsilon)\to 0$ as
$\varepsilon\to 0$. \label{thm1.1}
\end{theorem}

  Zongming Guo and Jiay Li \cite{Z.G} studied sequences of smooth
solutions of \eqref{eq1.1}  having uniformly bounded energy, they
proved the following result.

\begin{theorem}
Let ${u_{i}}$ be a sequence of smooth solutions of \eqref{eq1.1}
such that
\[
{\|u_{i}\|_{{H}^{1}(\Omega)}+\|u_{i}\|_{{L}^{2n/(n-2)}(\Omega)}}
\]
is bounded. Let $ u_{\infty}$ be the weak limit of $ {u_{i}}$ in
$\textbf{H}^{1}(\Omega)$ $\cap \textbf{L}^{2n/(n-2)}(\Omega)$.
Then $ u_{\infty}$ is smooth and satisfies  equation \eqref{eq1.1}
outside a closed singular subset $\Sigma$ of $\Omega$. Moreover,
there exists $r_0 >0$ and $\varepsilon_0 >0$ such that
\[
\Sigma = \cap_{0 < r < r_0 } \big\{ x \in \Omega: \liminf_{i \to
\infty}E_{u_i}(x,r) \geq \varepsilon _0  \big\}.
\]
\label{thm1.2}
\end{theorem}

We define the sequence of Radon measures
\[
{\eta}_i:= (\frac{1}{2} { \vert{\nabla u_{i}} \vert}^2 +
\frac{n-2}{2n} {\vert{u_i}\vert}^{2n/(n-2)}) \,dx
\]
Assumption that the sequence $( { \Vert {\nabla u_i}
\Vert}_{\mathbf H^{1}(\Omega)} +{\Vert{u_i}\Vert}_{\mathbf{L}
^{2n/(n-2)}(\Omega)})_i$ is bounded, and up to a
subsequences,  we can assume that
$\eta_i  \rightharpoonup \eta$
in the sense of measures as $i \to \infty$. Namely, for
any continuous function $\phi$ with compact support in $\Omega $
\[
\lim_{i \to \infty} \int_\Omega \phi \,d\eta_i  =
\int_ \Omega \phi \,d\eta.
\]
 Fatou's Lemma then implies that we can decompose
\[
\eta = (\frac{1}{2} { \vert{\nabla u} \vert}^2 +
\frac{n-2}{2n} {\vert{u}\vert}^{2n \over n-2}) \,dx+\nu
\]
where $\nu $ is a nonnegative Radon measure. Moreover, we prove
that $\nu$ satisfies the following lemma.

\begin{lemma}\label{lem1.2}
 Let $\delta>0$ such that $B_\delta\subset \Omega$. Then we
 have \begin{itemize}

\item[(i)] $\Sigma\subset spt(\nu)$

\item[(ii)] There exists a measurable, upper-semi-continuous function
$\Theta$ such that
\[
\nu (x)= \Theta (x)   \mathcal{H}^{0} \lfloor \Sigma ,\quad  \text{for }
x \in \Sigma.
\]
\end{itemize}
 Moreover, there exists some constants $c$ and
$C >0$ (only depending on $n$ and $\Omega$) such that
\[
c  \varepsilon_0 < \Theta (x) < C   \quad \mathcal{H}^{0}-\mbox{a.e. in }  \Sigma
\]
where ${H}^{0}\lfloor \Sigma$ is the restriction to $\Sigma$ of
the Hausdorff measure and $\Theta$ is a measurable function on
$\Sigma$.
\end{lemma}

The main question we would like to address in the present paper
concerns the multiplicity $\Theta$ of the defect measure which has
been defined above. More precisely, we have proved the following
theorem.

\begin{theorem}
Let $\nu$ be the defect measure of the sequence
$(|\nabla{u_{i}}|^{2}+|u_{i}|^{2n/(n-2)})dx$ defined above.
Then $\nu$ is quantized. That is,
 for a.e  $x\in \Sigma$,
\begin{equation}
\Theta(x)=\sum_{j=1}^{j=N_{x}}{{\|\nabla{v_{x,j}}\|^2_{{L}^{2}(\Omega)}+
\|{v_{x,j}}\|^{2n/(n-2)}_{{{L}^{2n/(n-2)}}(\Omega)}}}
\end{equation}
where $N_{x}$ is a positive integer and where the functions
$v_{x,j}$ are solutions of $\Delta v +  v^{\frac{n+2}{n-2}} =0$
which are defined on $\mathbb{R}^n$, issued from $(u_{i'})$ and
that concentrate at $x$ as $i\to \infty$. \label{thm1.3}
\end{theorem}

 The sentence ``issued from $(u_{i'})$ and that
concentrate at $x$ as $i \to \infty $''  means that there
are sequences of conformal maps ${\psi}_{j}^{i}$ , a finite family
of balls $(B^{l}_{i,j})_{l}$ such that the pulled back function
\[
{\tilde{u}_{i,j}} = ({\psi}_{j}^{i})^*    u_{i'}
\]
satisfies
\begin{gather*}
{\tilde{u}_{i,j}} \to v_j \quad \mbox{strongly in} \quad
\mathbf{L}^{2}({\mathbb R^n} \setminus {\cup}_{l}{B^{l}_{i,j})}),\\
\nabla{{\tilde{u}_{i,j}}} \to \nabla{v_j} \quad
\mbox{strongly in} \quad \mathbf{L}^{2}({\mathbb R^n} \setminus
{\cup}_{l} {B^{l}_{i,j})}  )
\end{gather*}

In the context of Yang-Mills fields in dimension $n \geq 4$ a
similar concentration result has been proven by  Rivi\`ere
\cite{T.R}. More precisely,  Rivi\`ere has shown that, if
$(A_i)_i$ is a sequence of Yang-Mills connections such that $(
\Vert {\nabla}_{A}{\nabla}_{A} F_{A}
\Vert_{\mathbf{L}^{1}(B_{1}^n) } )_i$ is bounded, then the
corresponding defect measure $\nu = \Theta  \mathcal{H}^{n-4}
\lfloor{\Sigma}$ of a sequence of smooth Yang-Mills connections is
quantized.

 The proof of Theorem \ref{thm1.3} uses technics introduced
by  Lin and  Rivi\`ere in their study of Ginzburg-Landau
vortices \cite{L.R} and also the technics developed by
Rivi\`ere in \cite{Z.G}. These technics use as an essential tool
the Lorentz spaces, more specifically the
$\mathbf{L}^{2,\infty}$-$\mathbf{L}^{2,1}$ duality \cite{L.T}.

\

 This paper is organized in the following way: In Section
2 we establish first a monotonicity formula for smooth solutions
of problem \eqref{eq1.1} which allows us to prove an
$\varepsilon$-regularity Theorem. Then, we prove Theorem
\ref{thm1.1} and Lemma \ref{lem1.2}. While Section 3 is devoted to
the proof of our main result, Theorem \ref{thm1.3}.

\section{A monotonicity Inequality}

In this section, we establish a monotonicity formula for smooth
solutions of problem \eqref{eq1.1}. Using Pohozaev identity:
Multiplying \eqref{eq1.1} by
$x_{i} \frac{\partial{u}}{\partial{x_{i}}}$ (summation over i is
understood) and integrating over B(x,r), the ball centered at $x$
of radius r, we obtain
\[
-\int_{B(x,r)}{x_{i}\frac{\partial{u}}{\partial{x_{i}}}}{\Delta{u}}\,dy
=-\int_{B(x,r)}{x_{i}\frac{\partial{u}}{\partial{x_{i}}}}u
|{u}|^{4/(n-2)}\,dy
\]
By Green formula, we get
\begin{equation}
\begin{aligned}
&\frac{n-2}{2} \int_{B(x,r)}{|u|^{2n/(n-2)}}dy
-\frac{n-2}{2}
 \int_{B(x,r)}{|\nabla{u}|^{2}}dy\\
&-\frac{n-2}{2n} \int_{\partial{B(x,r)}}{|u|^{2n/(n-2)}}\,ds
 + \frac{1}{2} r
\int_{\partial{B(x,r)}}{|\nabla{u}|^{2}}ds \\
&=  r \int_{\partial{B(x,r)}}{|\frac{\partial{u}}{\partial{r}}|^{2}}dy
\end{aligned}
\label{eq2.1}
\end{equation}
On the other hand, multiplying \eqref{eq1.1} by $u$ and
integrating over $B(x,r)$, we get
\begin{equation}
\int_{B(x,r)}{|\nabla{u}|^{2}}dy
-\int_{\partial{B(x,r)}}{u\frac{\partial{u}}{\partial{r}}}ds
=\int_{B(x,r)}{|u|^{2n/(n-2)}}dy \label{eq2.2}
\end{equation}
Deriving (\ref{eq2.2}) with respect to $r$, we obtain
\begin{equation}\int_{\partial{B(x,r)}}{|\nabla{u}|^{2}}dy
-\frac{d}{dr} \int_{\partial{B(x,r)}}{u\frac{\partial{u}}{\partial{r}}}ds
=\int_{\partial{B(x,r)}}{|u|^{2n/(n-2)}}dy \label{eq2.3}
\end{equation}
Combining \eqref{eq2.1}, \eqref{eq2.2} and \eqref{eq2.3}, we get
\begin{align*}
&-\frac{r}{n} \int_{\partial B(x,r)}|u|^{2n/(n-2)}\,ds \\
&= \frac{1}{2}r \frac{d}{dr}\int_{{\partial}B(x,r)}
u\frac{{\partial}u}{{\partial}r}\,ds - r\int_{{\partial}B(x,r)}
|\frac{{\partial}u}{{\partial}r}|^{2}\,dy +
r^{-1} u \frac{{\partial}u}{{\partial}r}\,ds  .
\end{align*}
Moreover, we have that
\begin{align*}
\frac{d^2}{dr^2} (\int_{{\partial}B(x,r)}u^2 \,ds )
&=\frac{d}{dr}(
2\int_{{\partial}B(x,r)}u\frac{{\partial}u}{{\partial}r}\,ds +
\frac{n-1}{r}\int_{{\partial}B(x,r)}u^2\,ds ) \\
 &=(n-1)\Big[\frac{2}{r}\int_{{\partial}B(x,r)}u\frac{{\partial}u}{{\partial}r}\,ds
+(\frac{n-1}{r^2}-\frac{1}{r^2})\int_{{\partial}B(x,r)}u^2\,ds
\Big]\\
&\quad +2\frac{d}{dr}\int_{{\partial}B(x,r)}u\frac{{\partial}u}{{\partial}r}\,ds\\
&=\frac{n-1}{r}\Big[2\int_{{\partial}B(x,r)}u\frac{{\partial}u}{{\partial}r}\,ds
+\frac{n-2}{r}\int_{{\partial}B(x,r)}u^2 \,ds\Big]\\
&\quad +
2\frac{d}{dr}\int_{{\partial}B(x,r)}u\frac{{\partial}u}{{\partial}r}\,ds.
\end{align*}
Hence
\begin{align*}
&\frac{1}{n} \frac{d}{dr}\int_{B(x,r)}|u|^{2n/(n-2)}\,dy
+\frac{1}{n} \frac{d^2}{dr^2} \int_{{\partial}B(x,r)}u^2
\,ds\\
 &=  \int_{{\partial}B(x,r)}(
|\frac{{\partial}u}{{\partial}r}|^2 +
\frac{2n-3}{2r} u \frac{{\partial}u}{{\partial}r} +
\frac{(n-1)(n-2)}{4}r^{-2}u^{2})\,ds.
\end{align*}
Moreover
\begin{align*}
&\frac{d}{dr}( \frac{1}{r}\int_{{\partial}B(x,r)}u^2 ds )\\
&=-\frac{1}{r^2}\int_{{\partial}B(x,r)}u^2 ds +
\frac{2}{r}\int_{{\partial}B(x,r)}u\frac{{\partial}u}{{\partial}r}ds+ \frac{n-1}{r^2}\int_{{\partial}B(x,r)}u^2 ds  \\
&= \frac{n-2}{r^2}\int_{{\partial}B(x,r)}u^2 ds +
\frac{2}{r}\int_{{\partial}B(x,r)}u\frac{{\partial}u}{{\partial}r}ds.
\end{align*}
We obtain
\begin{align*}
&\frac{d}{dr}\Big[
\frac{1}{n}\int_{B(x,r)}|u|^{2n/(n-2)}dy
+\frac{1}{n}\frac{d}{dr}\int_{{\partial}B(x,r)}u^2 ds
-\frac{1}{n}\frac{1}{r}\int_{{\partial}B(x,r)}u^2 ds \Big]\\
&=\int_{{\partial}B(x,r)}(|\frac{{\partial}u}{{\partial}r}|^2+
(n-2) r^{-1}u|\frac{{\partial}u}{{\partial}r}|
+\frac{(n-2)^2}{4}r^{-2}u^{2})ds\\
&=\int_{{\partial}B(x,r)}(\frac{{\partial}u}{{\partial}r}
 +\frac{n-2}{2}r^{-1}u)^2 ds \geq 0
\end{align*}
We  conclude that
\begin{equation}
E_u (x,r) = \frac{1}{n}\int_{B(x,r)}|u|^{2n/(n-2)}dy +
\frac{1}{n}\frac{d}{dr}\int_{B(x,r)}u^2 ds +
\frac{1}{n}r^{-1}\int_{B(x,r)}u^2 ds \label{eq2.4}
\end{equation}
is a nondecreasing function of $r$.
  Using the fact that
\[
\int_{B(x,r)}|u|^{2n/(n-2)}dy -
\int_{{\partial}B(x,r)}|\nabla u|^2 dy  = -
\int_{{\partial}B(x,r)} u\frac{{\partial}u}{{\partial}r}ds,
\]
one can easily get
 \begin{align*}
 E_u(x,r) &=
\frac{1}{n}\int_{B(x,r)}|u|^{2n/(n-2)}dy +
\frac{1}{4}\frac{d}{dr}\int_{{\partial}B(x,r)}u^2 ds -
\frac{1}{4}r^{-1}\int_{{\partial}B(x,r)}u^2 ds\\
&= \frac{\frac{n}{2}}{n}\int_{B(x,r)}|u|^{2n/(n-2)}dy
+\frac{1-\frac{n}{2}}{n}\int_{B(x,r)}|u|^{2n/(n-2)}dy\\
&\quad +\frac{1}{4}\frac{d}{dr}\int_{{\partial}B(x,r)}u^2 ds
-\frac{1}{4}r^{-1}\int_{{\partial}B(x,r)}u^2 ds \\
&= \frac{1}{2}\int_{B(x,r)}|{\nabla}u|^2 dy -
\frac{1}{2}\int_{{\partial}B(x,r)}
u\frac{{\partial}u}{{\partial}r}ds -
\frac{n-2}{2n}\int_{B(x,r)}|u|^{2n/(n-2)}dy
\\
&\quad +\frac{1}{4}\frac{d}{dr}\int_{{\partial}B(x,r)}u^2 ds
-\frac{1}{4}r^{-1}\int_{{\partial}B(x,r)}u^2 ds \\
&= \frac{1}{2}\int_{B(x,r)} (  |{\nabla}u|^2 -
\frac{n-2}{2n}|u|^{2n/(n-2)}  ) dy  + \frac{1}{4}
\frac{d}{dr}\int_{{\partial}B(x,r)}u^2 \,ds
\\
&\quad -\frac{1}{4}r^{-1}\int_{{\partial}B(x,r)}u^2 ds -
\frac{1}{2}\int_{{\partial}B(x,r)}
u\frac{{\partial}u}{{\partial}r} ds .
\end{align*}
We obtain an equivalent formulation of $ E_u(x,r)$
\begin{equation}
E_{u}(x,r)=\frac{1}{2}\int_{B(x,r)} ( |{\nabla}u|^2
-\frac{n-2}{2n}|u|^{2n/(n-2)} dy
+\frac{n-2}{4}r^{-1}\int_{{\partial}{B(x,r)}}u^2 ds
\label{eq2.5}
\end{equation}
Moreover, using the fact that
\[
\frac{d}{dr}\int_{\partial{B(x,r)}}{u^{2}}ds
=2\int_{\partial{B(x,r)}}{u\frac{\partial{u}}{\partial{r}}}ds
+\frac{n-1}{r}\int_{\partial{B(x,r)}}{u^2}
\] we obtain
\begin{align*}
\frac{1}{r}\int_{\partial{B(x,r)}}{u^2}\,ds
&=\frac{1}{n-1}\frac{d}{dr}\int_{\partial{B(x,r)}}{u^{2}}\,ds
-\frac{2}{n-1}\int_{\partial{B(x,r)}}{u\frac{\partial{u}}{\partial{r}}}\,ds\\
&=\frac{1}{n-1}\frac{d}{dr}\int_{\partial{B(x,r)}}{u^{2}}\,ds\\
&\quad +\frac{2}{n-1}\Big[\int_{B(x,r)}{|u|^{2n/(n-2)}}\,dy
-\int_{B(x,r)}{|\nabla{u}|^{2}}\,dy\Big]
\end{align*}
Then $E_{u}(x,r)$ can also be written
\begin{align*}
&E_{u}(x,r)\\
&=\frac{1}{2(n-1)}\int_{B(x,r)}(|\nabla{u}|^{2}
+\frac{n-2}{n}|u|^{2n/(n-2)})\,dy
+\frac{n-2}{4(n-1)}\frac{d}{dr}\int_{\partial{B(x,r)}}{u^{2}}\,ds.
\end{align*}

\begin{proof}[Proof of Lemma \ref{lem1.1}]
To prove that $(x,r) \mapsto  E_{u}(x,r)$ is continuous it
suffices to prove that
\[
(x,r) \mapsto\int_{\partial{B(x,r)}}{u^{2}}ds
\]
is continuous with respect to $x$ and $r$. We have
\[
\int_{\partial{B(x,r)}}{u{\frac{\partial{u}}{\partial{r}}}}ds =
\int_{B(x,r)}{|\nabla{u}|^{2}}-
\int_{B(x,r)}{|u|^{2n/(n-2)}}dy
\]
 Thus
 $(x,r)\mapsto{\int_{\partial{B(x,r)}}}{u{\frac{\partial{u}}{\partial{r}}}}$
is continuous, and this allows to get the conclusion.

Now, to prove that $E_u$ is positive, we proceed by contradiction.
If the result is not true, then there would exists $x \in \Omega $
and $R >  0 $  such that $E_{u}(x,R)< 0$. For almost every $y$
in some neighborhood of $x$, we
 have
\[
\lim_{r\to 0}{\int_{\partial{B(x,r)}}u{\frac{\partial{u}}{\partial{r}}}}\,ds
=0
\]
 integrating $E_{u}(x,r)$ over the interval
$[0,R]$ and using the fact that r $\mapsto{E_{u}(x,r)} $ is
increasing, we obtain
\begin{align*}
{\int_{0}}^{R}{E_{u}(y,r)}dr
&=\frac{1}{2(n-1)}{\int_{0}}^{R}{dr \int_{B(y,r)}}(|\nabla{u}|^{2}+\frac{n-2}{2n}{|u|^{2n/(n-2)}})dx\\
&\quad + \frac{n-2}{4(n-1)}\int_{\partial{B(y,R)}}{u^{2}}\,ds \\
&\leq {R  E_{u}(y,R)} < 0
\end{align*}
which is not possible. This proves Lemma \ref{lem1.1}.
\end{proof}

\begin{lemma}
There exist $r_{0}>0$ and some constant $c>0 $, depending only on
$n$, such that
\[
\int_{B(x,r)}(|\nabla{u}|^{2}+{|u|^{2n/(n-2)}})\,dy < c E_{u}(x,r)
\]
 for any r $<r_0/2$.
\label{lem2.1}
\end{lemma}

\begin{proof}
Using the fact that $(x,r) \mapsto E_{u}(x,r)$ is nondecreasing,
we have
 \begin{align*}
 r E_{u}(x,r)&\geq {\int_{0}^{r}} {E_{u}(x,s)} \,ds\\
&=\frac{1}{2n-2}\int_{0}^{r}ds\int_{B(x,s)}(|\nabla{u}|^{2}
  +\frac{n-2}{n}|u|^{2n/(n-2)})\,dy\\
&\quad +\frac{n-2}{4(n-1)}\int_{0}^{r}ds{\int_{\partial{B(x,s)}}}u^{2}
\,d\sigma\\
&\geq {\frac{1}{2(n-1)}\frac{n-2}{n}\int_{\frac{r}{2}}^{r}ds
\int_{B(x,s)}(|\nabla{u}|^{2}+|u|^{2n/(n-2)})dy}\\
&\geq C(n)\frac{r}{2}\int_{B(x,\frac{r}{2})}(|\nabla{u}|^{2}+|u|^{2n/(n-2)})
\,dy
\end{align*}
where $C(n)$ is a positive constant depending only on $n$. This
gives the desired result.
\end{proof}

As a consequence of Lemma \ref{lem2.1}, we have the following result.

\begin{lemma}
Assume that there exist $x_{0}$ and $r_{0} > 0 $ such that
$E_{u}(x_{0},r_{0})\leq \varepsilon $
then
\[
\int_{B(x,r)}(|\nabla{u}|^{2}+\frac{n-2}{n}|u|^{2n/(n-2)})\,dy\leq
C \varepsilon \quad \forall \quad 0<r<2r_{0}\] where $C$ is a
positive constant depending only on $n$. \label{lem2.2}
\end{lemma}

\begin{proof} Let $x_{0}$ and $r_{0}$ be such that
$E_{u}(x_{0},r_{0})\leq \varepsilon $ and let $0<r<r_{0}$, then for all
$ x \in B(x_{0},\frac{r}{2}) $ we have
\[
B(x,\frac{r}{2})\subset{B(x_{0},r)} \subset{B(x_{0},r_{0})}
\]
Thus
\begin{align*}
E_{u}(x_{0},r_0)
&\geq \frac{n-2}{2n(n-1)}\int_{B(x,\frac{r}{2})}|u|^{2n/(n-2)}\,dy\\
&\quad + \frac{1}{2(n-1)}\int_{B(x,\frac{r}{2})}|\nabla{u}|^{2}dy
 +\frac{n-2}{4(n-1)}\frac{d}{dr}\int_{\partial{B(x_{0},r)}}u^{2}\,ds\\
&\geq
\frac{1}{2(n-1)}\int_{B(x,\frac{r}{2})}\left(|u|^{2n/(n-2)}+\nabla{u}|^{2}\right)
 dy + \frac{n-2}{4(n-1)}\frac{d}{dr}\int_{\partial{B(x_{0},r)}}u^{2}\,ds
\end{align*}
Integrating between $0$ and $r$, we obtain
\begin{align*}
&r E_{u}(x_{0},r_{0})\\
&\geq \frac{1}{2(n-1)} \int_{0}^{r}ds \int_{B(x,
 \frac{s}{2})}(|u|^{2n/(n-2)}+\nabla{u}|^{2})\,dy+\frac{n-2}{4(n-1)}\int_{\partial{B(x_{0},r)}}u^{2}\,ds\\
&\geq \frac{1}{2(n-1)}\int_{0}^{r}ds\int_{B(x,\frac{s}{2})}
 (|\nabla{u}|^{2}+|u|^{2n/(n-2)} )\,dy\\
 &\geq \frac{1}{2(n-1)}\int_{\frac{r}{2}}^{r}\,ds
   \int_{B(x,\frac{s}{2})}(|\nabla{u}|^{2}+|u|^{2n/(n-2)})\,dy\\
&\geq \frac{1}{2(n-1)}\frac{r}{2}\int_{B(x,\frac{r}{2})}
 (|\nabla{u}|^{2}+|u|^{2n/(n-2)})\,dy.
 \end{align*}
 Then
 \[
 E_{u}(x_{0},r_{0})\geq\frac{1}{4(n-1)}
 \int_{B(x,\frac{r}{2})}(|\nabla{u}|^{2}+|u|^{2n/(n-2)} )\,dy,
 \]
 thus
\[
\int_{B(x,r)}(|\nabla{u}|^{2}+|u|^{2n/(n-2)})\,dy \leq{C
\varepsilon}      \quad \forall  r<2r_{0}.
\]
This proves the desired result.
\end{proof}

\begin{proof}[Proof of Theorem \ref{thm1.1}]

Without loss of generality, we can assume that $x_{0}=0$ and we
denote by $B_{r_{0}}$ the ball of radius $r_{0}$ centered at
$x_{0}=0$ .

We use the idea of Schoen \cite{Sch}. For $r<r_{0}$, we define
\[
F(y)=(\frac{r}{2}-|y|)^{(n-2)/2}u(y)
\]
Clearly $F$ is continuous over $B_{\frac{r}{2}}$, then there exist
$y_{0}\in B_{\frac{r}{2}}$ such that
\[
F(y_{0})=\max_{y \in{B_{\frac{r}{2}}}}{(\frac{r}{2}
 -|y|)^{(n-2)/2}u(y)}=(\frac{r}{2}-|y_{0}|)^{(n-2)/2}u(y_{0})
\]
Let $0<\sigma<\frac{r}{2}$, for all $y\in B_{\sigma}$, we have
\[
u(y)\leq\frac{(\frac{r}{2}-|y_{0}|)^{(n-2)/2}}{(\frac{r}{2}
-|y|)^{(n-2)/2}} u(y_{0})
\]
 Then
\[
\sup_{y\in {B_{\sigma}}}{u(y)}\leq\frac{(\frac{r}{2}
-|y_{0}|)^{(n-2)/2}}{(\frac{r}{2}-|y|)^{(n-2)/2}}\sup_{y\in
{B_{\sigma_{0}}}}u(y)
\]
where $\sigma_{0}=|y_{0}|$. Let $y_{1}\in
B_{\sigma_{0}}$ be such that
\[
u(y_{1})=\sup_{y\in {B_{\sigma_{0}}}}u(y)
\]
We claim that
\[
u(y_{1})\leq\frac{2^{(n-2)/2}}{(\frac{r}{2}-|y_{0}|)^{(n-2)/2}}.
\]
Indeed, on the contrary case, we get
\[
(u(y_{1}))^{-2/(n-2)}\leq\frac{1}{2}(\frac{r}{2}-|y_{0}|)
\]
Let $\mu=(u(y_{1}))^{-2/(n-2)}$. We have
\[
B_{\mu}(y_{1})\subset{B_{\frac{\sigma_{0}+\frac{r}{2}}{2}}}
\]
($|z-y_{1}|<\mu $ take $|z|<\frac{\frac{r}{2}+|y_{0}|}{2}$).
Hence
\[
\sup_{y\in
{B_{\mu}}(y_{1})}u(y)\leq\frac{(\frac{r}{2}-|y_{0}|)^{(n-2)/2}}
{(\frac{\frac{r}{2}-|y_{0}|}{2})^{(n-2)/2}}u(y_{1})
=2^{(n-2)/2}u(y_{1})
\]
Let
$ v(x)=\mu^{(n-2)/2}u(\mu{x}+y_{1})$.
Easy computations shows that $v$ satisfies
\begin{align*}
\Delta{v^{2n/(n-2)}}
&= \frac{2n}{n-2}\Big[\frac{n+2}{n-2}v^{4/(n-2)}|\nabla{v}|^{2}
 +v^{\frac{n+2}{n-2}}\triangle{v}\Big]\\
&\geq \frac{2n}{n-2}v^{\frac{n+2}{n-2}}\triangle{v}=
-\frac{2n}{n-2}v^{2\frac{n+2}{n-2}}
\end{align*}
 On the other hand
 \[
v^{2n/(n-2)}(0)
 =  {\mu}^{\frac{n-2}{2}\frac{2n}{n-2}}u^\frac{2n}{n-2}(y_{1})=1 .
 \]
Moreover, we have
\begin{align*}
\sup_{B_{1}}v(x)&={\mu}^{(n-2)/2}\sup_{B_{1}}u(\mu{x}+y_{1})\\
&= {\mu}^{(n-2)/2}\sup_{B_{\mu}(y_{1})}u(x)\\
&\leq {\mu}^{(n-2)/2}{2}^{(n-2)/2}u(y_{1})={2}^{(n-2)/2}.
\end{align*}
Then
$\sup_{B_{1}}v^{2n/(n-2)}\leq{2}^n$.
Therefore,
\[
-\Delta{v^{2n/(n-2)}}\leq{C(n)v^{2n/(n-2)}} .
\]
We conclude that
\[
1=v^{2n/(n-2)}(0)\leq{C\int_{B_{1}}v^{2n/(n-2)}(x)dx}
= C{\mu^n}\int_{B_{\mu}}{u^{2n/(n-2)}}(x)dx\leq C\varepsilon .
\]
For $\epsilon$ sufficiently small, we derive a contradiction.
It follows that
\[
\sup_{B_{\frac{r}{2}}}u(y)
\leq{\frac{(\frac{r}{2}-|y_{0}|)^{(n-2)/2}}{(\frac{r}{2}-|y|)^{(n-2)/2}}
\cdot
\frac{2^{(n-2)/2}}{(\frac{r}{2}-|y_{0}|)^{(n-2)/2}}}=\frac{2^{(n-2)/2}}{(\frac{r}{2}-|y|)^{(n-2)/2}}.
\]
For $|y|<r/4$, we have
\[
\sup_{B_{\frac{r}{4}}}u(y)\leq C(n)/ r^{(n-2)/2}
\]
This in turns proves the Theorem \ref{thm1.2}.
\end{proof}

\begin{proof}[Proof of Lemma \ref{lem1.2}]

We keep the above notations. To show (i), suppose $x_0 \in B_1
\setminus \Sigma$, then there exists $r_1 > 0$ such that
\[
\liminf_{i \to \infty}  E_{u_i}(x_0,r_1) <  \varepsilon_0.
\]
Then, we may find a sequence $n_j\to \infty$ as $j \to \infty$ such
that
\[
\sup_{n_j}E_{u_{n_j}}(x_0,r_1) < \varepsilon_0.
\]
We deduce from the $\varepsilon$-regularity Theorem (Theorem
\ref{thm1.1}) that
\[
\sup_{n_j}\sup_{x\in B_{\frac{r_1}{16}}(x_0)}|u_{n_j}|\leq
\frac{C}{r_1^{(n-2)/2}}.
\]
for some constant $C$ depending only on $n$. Then
\[
u_{n_j}\to u \quad\text{in  } C^1(B_{\frac{r_1}{16}}(x_0))
\]
a similar argument allows to show that
\[
\nabla u_{n_j}\to \nabla u \quad\text{in }
 C^1(B_{\frac{r_1}{16}}(x_0))
\]
Then
\[
\mu_{n_j}:=\left(\frac{1}{2} |\nabla u_{n_j}|^2+ \frac{n-2}{2n}
u_{n_j}^{2n/(n-2)}\right)\,dx\to \left(\frac{1}{2} |\nabla
u|^2+\frac{n-2}{2n} u^{2n/(n-2)}\right)\,dx
\]
as radon measure. Hence $\nu = 0$ on $B_{\frac{r_1}{16}}(x_0)$ i.e
$x_0 \notin \mathop{\rm supp}(\nu)$ and then we deduce that
$\mathop{\rm supp}(\nu)\subset \Sigma$.

To show (ii), let us first recall some properties of the function
$E_u(x,r)$ that has been defined above:

\noindent $\bullet$ For all $x \in \Omega$, there exists $r_0 > 0$ and a
constant $ C>0$ such that
\[\int_{B(x,r)}(\frac{1}{2} |\nabla u|^2+ \frac{n-2}{2n}|u|^{2n/(n-2)})
< C  E_u(x,r_0) \quad \forall   r<\frac{r_0}{2}.
\]
This is explained in the proof of Lemma \ref{lem1.1}.

\noindent$\bullet$ Using the fact that $E_u(x,.)$ is increasing on $r$
together with the fact that
\[
\lim_{r\searrow 0} E_u(x,r)=0 \quad \mathcal{H}^{0}-\text{a.e. }  x\in\Omega
\]
we deduce that for $\mathcal{H}^{0}$-a.e. $x\in\Sigma$,
$\lim_{r \searrow 0} \int_{B(x,r)} \nu$  exists.
and  the density $\Theta (\eta, .)$ defined by
\begin{equation}
\Theta (\eta ,x):= \lim_{r\searrow 0} \eta(B_r(x))
\end{equation}
exists for every $x\in \Omega$. Moreover, for $\mathcal{H}^{0}$-a.e.
$x\in\Omega$, $\Theta_u(x)=0$, where
\begin{equation}
\Theta_u(x):= \lim_{r \searrow 0}
\int_{B(x,r)}(\frac{1}{2} {\vert \nabla u \vert}^2
+ \frac{n-2}{2n} {\vert u \vert}^{2n \over {n-2}}) \,dy.
\end{equation}
 Now, for $r$
sufficiently small and $i$ sufficiently large
\begin{equation}
\int_{B(x,r)} \frac{1}{2}|\nabla u_i
|^2+\frac{n-2}{2n} u_i^{2n/(n-2)} \leq
 C E_{u_i}(x,r) \leq   C(\Lambda,\Omega) \label{eq4.13}
\end{equation}
where $\Lambda$ is given above and $ C(\Lambda,\Omega)$
is a constant depending only on $\Lambda $ and $\Omega$. Hence
\begin{equation}
 \eta(B(x,r)) \leq  C(\Lambda,\Omega) \quad\text{for } x
\in B_1^n \label{eq4.14}
\end{equation}
In particular, this implies that $\eta \lfloor \Sigma $ is
absolutely continuous with respect to $\mathcal{H}^{0}\lfloor \Sigma $.
Applying Radon-Nikodym's Theorem \cite{E.G}, we conclude that
\begin{equation}
\eta \lfloor \Sigma = \Theta (x) \mathcal{H}^{0} \lfloor \Sigma \quad
\text{for } \mathcal{H}^{0}\text{-a.e. } x \in \Sigma
\end{equation}
Using \ref{eq4.13} we conclude that
 \begin{equation}
\nu (x) = \Theta (x) \mathcal{H}^{0}\lfloor \Sigma
\end{equation}
for a $\mathcal{H}^{0}$-a.e. $x \in \Sigma$
(recall that $\eta =(\frac{1}{2}{\vert \nabla u \vert }^2 +\frac{n-2}{2n}
{\vert u \vert}^{2n \over {n-2}})\,dx + \nu $ and
$\mathop{\rm supp}( \nu) \subset \Sigma $). The  estimate on
$\Theta$ follows from \ref{eq4.14}.
\end{proof}

For any $ y\in B^n_1$ and any sufficiently small $\lambda > 0$, we
define the scaled measure $\eta_{y,\lambda}$ by
\begin{equation}
\eta_{y,\lambda}(x) := \eta (y + \lambda  x )
\end{equation}
We have the following lemma.

\begin{lemma} \label{lem4.1}
Assume that $(\lambda_j)_j$ satisfies $\lim_{j\to \infty}\lambda_j=0$.
Then, there exist a subsequence
$({{\lambda}_{j'}} )_{j'}$ and a Radon measure $\chi$ defined on
$\Omega$, such that
${\eta_{y,\lambda_{j'}}} \rightharpoonup \chi $
in the sense of measures.
\end{lemma}

\begin{proof}
 For each $i \in {\mathbb N}$, we define
the scaled function $u_{i,y,\lambda}$ by
\begin{equation}
u_{i,y,\lambda}(x):={\lambda}^{{n-2 \over 2}} u_i(\lambda x +
y)\quad\text{for } y \in  B^n_1.
\end{equation}
Then $u_{i,y,\lambda}$ is a solution of
\[
-\Delta{u}= u|u|^{4/(n-2)} \quad \text{on } B^n_1.
\]
 In addition, for any  $r > 0$ sufficiently small, we have
\begin{equation}
\begin{aligned}
&\int_{B_r(0)}\left(\frac{1}{2}{\vert \nabla u_{i,y,\lambda}
\vert}^2 +\frac{n-2}{2n} {\vert u_{i,y,\lambda}\vert}^{{2k
\over {k-2}}}\right) \,dx \\
&=  \int_{B_{\lambda r}(y)} \left(\frac{1}{2} {\vert \nabla u_i
\vert}^2 + \frac{n-2}{2n}  {\vert u_i\vert}^{{2n \over
{n-2}}}\right) \,dx \leq C(\Lambda , \Omega).
\end{aligned}
\label{eq419}
\end{equation}
Finally for fixed $\lambda$,
\begin{align*}
&\left(\frac{1}{2} {\vert \nabla u_{i,y,\lambda} \vert}^2  +
 \frac{n-2}{2n} {\vert u_{i,y,\lambda}\vert}^{2n/(n-2)} \right) (x) \,dx
\\
&={\lambda}^n  \left(\frac{1}{2} {\vert \nabla u_i \vert}^2
 - \frac{n-2}{2n}  {\vert u_i \vert}^{2n/(n-2)}\right)(\lambda  x +y)\,dx
\\
& \rightharpoonup \eta(\lambda x +y) = \eta_{y,\lambda}(x)
\end{align*}
in the sense of measures as $i\to \infty$. On the other
hand letting $i$ tends to infinity in (\ref{eq419}), we conclude
that for any $r>0$
\begin{equation}
\eta_{y , \lambda} (B_r(0)) \leq  C(\Omega,\Lambda).
\end{equation}
Hence, we may find a subsequence $\{\lambda'_j \}$ of $\{
\lambda_j\}$ and a Radon measure $\chi$ such that $\eta_{y ,
\lambda'_j}$ converge weakly to $\chi$ as Radon measure on
$\Omega$. Then
\[
\lim_{j\to \infty} \lim_{i\to \infty} \left(  \frac{1}{2} {\vert
\nabla u_{i,y,\lambda'_j} \vert}^2 +\frac{n-2}{2n} {\vert
u_{i,y,\lambda'_j} \vert}^{2n \over {n-2}}\right)\,dx = \lim_{j\to
\infty} \eta_{y,\lambda'_j}(x) = \chi
\]
Using a  diagonal subsequence argument, we may find a subsequence
$i_j \to\infty$, such that
\[
\lim_{j\to \infty} \left(\frac{1}{2}  {\vert \nabla
u_{i_j,y,\lambda'_j} \vert}^2 +\frac{n-2}{2n}  {\vert
u_{i_j,y,\lambda'_j} \vert}^{2n \over {n-2}}\right)dx  = \chi
\]
This proves the Lemma.
\end{proof}

\begin{remark}  \label{rmk1}\rm
Observe that
\[
\chi(B_r(0))= \lim_{j\to \infty } \eta_{y ,
\lambda'_j}(B_r(0))=\lim_{j\to \infty }\eta (B_{\lambda'_j
r}(y))=\Theta(\eta, y)
\]
In particular, we  deduce that $\chi(B_r(0))$ is independent of r.
\end{remark}

\section{Proof of Theorem \ref{thm1.3}}

The idea of the proof comes from Rivi\`ere \cite{T.R} in the
context of Yang-Mills Fields. To simplify notation and since the
result is local, we assume that $\Omega$ is the unit ball $B^n$ of
$\mathbb{R}^n$. Let $(u_{k}) $ be a sequence of smooth solutions
of \eqref{eq1.1} such that
\[
\Big( \|u_{k}\|_{\textbf{H}^{1}(\Omega)}
+\|{u_{k}}\|_{\textbf{L}^{2n/(n-2)}(\Omega)} \Big)
\]
 is bounded and let $\nu$ be the defect measure defined above.
We claim that  for $\delta>0$, we have
\begin{equation}
\lim_{k\to\infty} \sup_{y\in{B_{1}}(x_{0})}\int_{B_\delta(y_{0})}
\left(|u_{k}|^{2n/(n-2)}+|\nabla{u_{k}}|^{2}\right)\geq{\varepsilon(n)} \label{eq3.1}
\end{equation}
where $\varepsilon (n)$ is given by Theorem \ref{thm1.3}. Indeed if
(\ref{eq3.1}) would not hold, we have for $\delta>0$ and
$k\in\mathbb{N}$ large enough
\[
\sup_{y\in{B_{1}}(x_{0})}\int_{B_\delta(y_{0})}\left(|u_{k}|^{2n/(n-2)}+|\nabla{u_{k}}|^{2}\right)\leq{\varepsilon(n)}
\]
and by Theorem \ref{thm1.1} we have
\[
\|\nabla{u_{k}}\|_{\textbf{L}^{\infty}(B_\frac{\delta}{2}(y))}
\leq C(\epsilon)/ r^{n/2}
\]
This contradict the concentration phenomenon and the claim is
proved. We then conclude that there exists sequences
$\delta_{k}\to 0$ as $k\to\infty$ and $
(y_{k})\subset{B_{1}(x_{0})}$ such that
\begin{equation}
\begin{aligned}
\int_{B_{\delta_{k}}(y_{0})}\left(|u_{k}|^{2n/(n-2)}+|\nabla{u_{k}}|^{2}\right)dx
&=
\sup_{y\in{B_{1}}(x_{0})}\int_{B_{\delta_{k}(y_{0})}}\left(|u_{k}|^{2n/(n-2)}
+|\nabla{u_{k}}|^{2}\right)\,dx\\
&= \frac{\varepsilon(n)}{2}. \label{eq3.2}
\end{aligned}
\end{equation}
In other words, $y_{k}$ is located at a bubble of characteristic
size $\delta_{k}$. More precisely, if one introduces the function
\[
{\widetilde{u}}_{k}(x)=\delta_{k}^{(n-2)/2}
u_{k}(\delta_{k}x+y_{k});
\]
we have, up to a subsequence, that
\begin{gather*}
{\widetilde{u}}_{k} \to {u_{\infty}}\quad \text{in
}{\textbf{C}}^{\infty}_{\rm loc}(\mathbb{R}^n) \quad
\text{as }k \to \infty ,\\
\nabla{{\widetilde{u}}_{k}}\to{\nabla{u_{\infty}}} \quad\text{in }
{\textbf{C}}^{\infty}_{\rm loc}(\mathbb{R}^n)\quad  \text{as }  k
\to \infty\,.
\end{gather*}
Therefore,
\[-\Delta{u_{\infty}}=u_{\infty} |{u_{\infty}}|^{4/(n-2)} \quad
\text{in  } \mathbb{R}^n.
\]
This is the first bubble we detect. On the other hand, we have
clearly that
\begin{equation}
\int_{{\mathbb{R}}^n}\left(|u_{\infty}|^{2n/(n-2)}+|\nabla{u_{\infty}}|^{2}\right)\,dx
= \lim_{R\to\infty}
\lim_{k\to\infty}\int_{B_{R{\delta_{k}}}(y_{k})}
\left(|u_{k}|^{2n/(n-2)}+|\nabla{u_{k}}|^{2}\right)\,dx.
\label{eq3.3}
\end{equation}
Indeed:
\begin{align*}
&\lim_{R\to\infty} \lim_{k\to\infty}\int_{B_{R{\delta_{k}}}(y_{k})}
\left(|u_{k}|^{2n/(n-2)}+|\nabla{u_{k}}|^{2}\right) \,dx \\
&= \lim_{R\to\infty} \lim_{k\to\infty}\int_{B_{R}(0)}
\left(|u_{k}|^{2n/(n-2)}+|\nabla({u_{k})}|^{2}\right)(\delta_{k}x+y_{k})\,
\delta_{k}^n\,dx \\
&= \lim_{R\to\infty} \lim_{k\to\infty}\int_{B_{R}(0)}
\left(|{\delta_{k}}^{\frac{2-n}{2}}{\widetilde{u}_{k}(x)}|^{2n/(n-2)}+
|{\delta_{k}}^{\frac{2-n}{2}}{\delta_{k}}^{-1}\nabla{\widetilde{u}_{k}(x)}|^{2}\right)
 \delta_{k}^n\,dx \\
 &=
\lim_{R\to\infty} \lim_{k\to\infty}\int_{B_{R}(0)}
\left(|{\widetilde{u}_{k}(x)}|^{2n/(n-2)}+|\nabla{\widetilde{u}_{k}(x)}|^{2}\right)\,dx \\
&=
\lim_{R\to\infty}\int_{B_{R}(0)}\left(|u_{\infty}(x)|^{2n/(n-2)}+|\nabla{u}_{\infty}(x)|^{2}\right)\,dx \\
&=\int_{{\mathbb{R}}^n}\left(|u_{\infty}(x)|^{2n/(n-2)}
+|\nabla{u}_{\infty}(x)|^{2}\right)\,dx\,.
\end{align*}

Assume first that we have only one bubble of characteristic
$\delta_{k}$. We have shown that
\begin{equation}
\Theta=\lim_{k\to\infty}\int_{B_{1}^n(0)}\left(|\nabla{u_{k}}|^{2}
+|u_{k}|^{2n/(n-2)}\right)\,dx
=\int_{\mathbb{R}^n}\left(|\nabla{u_{\infty}}|^{2}+|u_{\infty}|^{2n/(n-2)}\right)\,dx,
\label{eq5.24}
\end{equation}
where $\Theta$ is defined above. It suffices to prove that
\begin{equation}
\lim_{R\to\infty} \lim_{k\to\infty}\int_{B_{1}^n(0)
\setminus{B_{R\delta_{k}(y_{k})}}}
\left(|{u}_{k}(x)|^{2n/(n-2)}+|\nabla{{u}_{k}(x)}|^{2}\right) \,dx
= 0\,. \label{eq3.4}
\end{equation}
In other words there is no ``neck'' of energy which is quantized.

To simplify notation, we assume that $y_{k} = 0$. We claim that
for any $\varepsilon > 0$ small enough, there exists  $R> 0$ and $k_{0}\in
\mathbb{N}$ such that for any $k\geq{k_{0}}$ and
$R\delta_{k}\leq{r}\leq{\frac{1}{2}}$, we have
\begin{equation}
\int_{B_{2r}^n(0)\setminus{B_{r(0)}}}\left(|{u}_{k}(x)|^{2n/(n-2)}
+|\nabla{{u}_{k}(x)}|^{2}\right)\,dx \leq\varepsilon \label{eq3.5}
\end{equation}
Indeed, if is not the case, we may find $\varepsilon_{0}>0$,  a subsequence
$k'\to\infty$ (Still denoted $k$ ) and a sequence $r_{k}$
such that
\begin{equation}
\begin{gathered}
\int_{B_{2r}^n(0)\setminus{B_{r}(0)}}\left(|{u}_{k}(x)|^{2n/(n-2)}
+|\nabla{{u}_{k}(x)}|^{2}\right)\,dx >\varepsilon_0, \\
\frac{r_k}{\delta_k}\to\infty \quad  \text{as }\quad k\to\infty
\end{gathered}
\label{eq3.6}
\end{equation}
Let $\alpha_k\to0$ such that $r_k/\alpha_k=o(1)$
and $\alpha_k r_k/\delta_k \to\infty$ and let
\[
v_k(x)=r_k^{(n-2)/2} u_k(r_k x)
\]
clearly  $v_k $ satisfies
\[
-\Delta{v_k}=v_k  |v_k|^{4/(n-2)} \quad\text{in }  B_{2\alpha_k}
\setminus{B_{\alpha_k}}
\]
Therefore,
\[
\int_{B_{2}^n(0)\setminus{B_{1}(0)}}\left(|{v}_{k}(x)|^{2n/(n-2)}
+|\nabla{{v}_{k}(x)}|^{2}\right)\,dx >\varepsilon(n)
\]
and then we have a second bubble. This contradict our assumption.

We deduce from (\ref{eq3.6}) and Theorem \ref{thm1.1} that for any
$\varepsilon <\varepsilon(n)$, there exist $R>0$ and  $ k_0 \in \mathbb{N}$ such that
for all $k\geq{k_0}$ and $|x|\geq{R\delta_k}$
\[
|\nabla{u_k}|(x)\leq C(\epsilon) /|x|^{n/2}
\]
where $C(\varepsilon)\to 0$ as $\varepsilon\to0$. Then
\begin{equation}
|\nabla{u_k}|^{2}(x)\leq C(\varepsilon)/|x|^n . \label{eq3.7}
\end{equation}
We define $E_\lambda^k$ by
 \[
E_\lambda^k = \mathop{\rm meas} \left\{x \in{\mathbb{R}^n} :
|\nabla{u_k}|(x)\geq\lambda\right\}
\]
We have
$ E_\lambda^k \leq C(\varepsilon)/ \lambda ^2$;
indeed
\[
\left\{x\in{\mathbb{R}^n} :
|\nabla{u}_k|(x)\geq\lambda\right\}\subset\{x\in{\mathbb{R}^n} :
|x|^n\leq{\frac{C(\varepsilon)}{{\lambda}^{2}}}\}
\]
and
\[
\mathop{\rm meas} \left\{x\in{\mathbb{R}^n} :
|x|^n\leq{\frac{C(\varepsilon)}{{\lambda}^{2}}}\right\}
\leq{\frac{C(\varepsilon)}{{\lambda}^2}}
\]
We deduce from (\ref{eq3.7}) that
 \begin{equation}
 \|\nabla{u_k}\|_{{\textbf{L}}^{2,\infty}(C_{B_{R{\delta}_k}})}\leq C(\varepsilon)
 \label{eq3.8}
 \end{equation}
where ${{\textbf{L}}^{2,\infty}}$ is the Lorentz space defined in
\cite{L.T}, the  weak ${\textbf{L}}^{2}$ space, and
$\|\cdot\|_{{\textbf{L}}^{2,\infty}}$ is the weak norm defined by
\[
\|f\|_{{\textbf{L}}^{2,\infty}}= \sup_{0< t <\infty} {t^{1/2}f^*(t)}
\]
where $f^*$ is the nonincreasing rearrangement of $|f|$. Indeed
\[
\|\nabla{u_k}\|_{{\textbf{L}}^{2,\infty}(C_{B_{R{\delta}_k}})} =
\sup_{0< t <\infty} {t^{1/2}(\nabla{u_k})^*(t)}
\]
 by definition,
\[
(\nabla{u_k})^*(t)=inf\{\lambda>0   /  E^k_\lambda   \leq t  \}
\]
For all $t >0$ such that $\frac{C(\varepsilon)}{\lambda^2}\leq t$, we have
$E^k_\lambda \leq t$. Then
\begin{align*}
\inf \left\{\lambda>0: E^k_\lambda \leq t  \right\}
&\leq \inf\left\{\lambda>0 : \frac{C(\varepsilon)}{\lambda^2} \leq t \right\}\\
&\leq \inf\left\{\lambda>0 : \lambda\geq \frac{(C(\varepsilon))^{1/2}}{t^{1/2}} \right\}\\
&=\frac{(C(\varepsilon))^{1/2}}{t^{1/2}}
\end{align*}
Hence $t^{1/2}(\nabla{u_k})^*(t) \leq C(\varepsilon)$
 and so
\begin{equation}
\|\nabla{u_k}\|_{{\textbf{L}}^{2,\infty}(C_{B_{R{\delta}_k}})}\leq
C(\varepsilon) \label{eq4.19}
\end{equation}
We claim that the sequence  $(\nabla{u_k})$ is uniformly bounded
in the Lorentz space ${\textbf{L}}^{2,1}(B^n_1)$ (see \cite{L.T}
for the definition). We prove this claim using an iteration
proceeding; Indeed, the sequence $(u_k)$ is bounded in
${\textbf{L}}^\frac{2n}{n-2}(B_1^n)$. Then
\[
\Delta{u_k}= - u_k |u_k|^{4/(n-2)}
\]
is bounded in ${\textbf{L}}^\frac{2n}{n+2}(B_1^n)$ which implies
by the elliptic regularity Theorem that the sequence $(u_k)$ is
bounded in ${\textbf{W}}^{2,\frac{2n}{n+2}}(B_1^n)$. Using the
imbedding Theorem for Sobolev spaces
\[
{\textbf{W}}^{m,p}(B_1^n)\subset {\textbf{W}}^{r,s}(B_1^n)  \quad
\text{if }  m \geq   r,\;   p \geq s \text{ and }
  m-\frac{n}{p}=r-\frac{n}{s}.
  \]
In particular, ${\textbf{W}}^{2,\frac{2n}{n+2}}(B_1^n)$   is
continuously imbedded in  ${\textbf{W}}^{1,2}(B_1^n)$. On the other
hand by Proposition 4 in \cite{L.T}, we have
\[
{\textbf{W}}^{1,2}(B_1^n)\hookrightarrow
{\textbf{L}}^{2^*,2}(B_1^n) =
{\textbf{L}}^{\frac{2n}{n-2},2}(B_1^n)
\]
continuously. We then deduce that
\[
\Delta{u_k}= - u_k |u_k|^{4/(n-2)}
\]
is bounded in ${\textbf{L}}^{\frac{2n}{n+2},\frac{2(n-2)}{n+2}}(B_1^n)$.
Here, we have used   the following lemma.

\begin{lemma}
If $f \in{\textbf{L}}^{p,q}(B_1^n)$ and
$\alpha\in\mathbb{Q^+}$, then
$f^\alpha\in{\textbf{L}}^{\frac{p}{\alpha},\frac{q}{\alpha}}(B_1^n)$.
\label{lem5.1}
\end{lemma}

\begin{proof}
In the case where $\alpha\in\mathbb{N}$, the
result follows from the fact that
\[
f \in{\textbf{L}}^{a,b}(B_1^n) \text{ and }
g \in{\textbf{L}}^{c,d}(B_1^n)\Rightarrow    f.g
\in{\textbf{L}}^{q,r}(B_1^n),
\]
 where $\frac{1}{q}=\frac{1}{a}+\frac{1}{b}$ and
$\frac{1}{r}=\frac{1}{c}+\frac{1}{b}$  (see \cite{B.W}). The
general case is a consequence of the fact that the increasing
rearrangement of the function $|f|^\beta$
 is equal to the puissance $\beta$ of the increasing rearrangement
of $|f|$ since $(f^\beta)^*$ is the only one function verifying
\[
\mathop{\rm meas} \{x\in{\mathbb{R}^n}  : f^{\beta}(x) \geq\lambda  \} =
\mathop{\rm meas} \{t>0  :  (f^{\beta})^*(x) \geq\lambda  \}
\]
 This in turns proves  Lemma \ref{lem5.1}.
\end{proof}

 Now, using  in \cite[Theorem 8]{L.T}, we deduce from (\ref{eq3.6})
that $(\nabla{u_k})$ is uniformly bounded in the space
${\textbf{L}}^{(\frac{2n}{n+2})^*,\frac{2(n-2)}{n+2}}(B_1^n)$ =
 ${\textbf{L}}^{2,\frac{2(n-2)}{n+2}}(B_1^n)$. Hence $(u_k)$ is bounded in
 ${\textbf{L}}^{2^*,\frac{2(n-2)}{n+2}}(B_1^n)$. Then
\[
\Delta{u_k}= - u_k |u_k|^{4/(n-2)}
\]
 is bounded in ${\textbf{L}}^{\frac{2n}{n+2},\frac{2(n-2)^2}{(n+2)^2}}(B_1^n)$.
Hence, again by \cite[Theorem 8]{L.T},  the sequence
$(\nabla{u_k})$ is bounded in
${\textbf{L}}^{2,\frac{2(n-2)^2}{(n+2)^2}}(B_1^n)$ and by elliptic
regularity Theorem
\[
\Delta{u_k}= - u_k |u_k|^{4/(n-2)}
\]
is bounded in
${\textbf{L}}^{\frac{2n}{n+2},\frac{2(n-2)^3}{(n+2)^3}}(B_1^n)$.
We obtain after $p$ iterations that
\[
\Delta{u_k}= - u_k |u_k|^{4/(n-2)}
\]
 is bounded in ${\textbf{L}}^{\frac{2n}{n+2},\frac{2(n-2)^p}{(n+2)^p}}(B_1^n)$.
We choose $p >0$ such that $6p >n$, we have in particular
$\frac{2(n-2)^p}{(n+2)^p}<1 $  which gives
\[
\Delta{u_k}= - u_k |u_k|^{4/(n-2)}
\]
 is bounded in
${\textbf{L}}^{\frac{2n}{n+2},1}(B_1^n)$. Here we have used the
fact that
\[
{\textbf{L}}^{p, q_1}(B_1^n)\subset{\textbf{L}}^{p, q_2}(B_1^n)
\quad\text{if }q_1 < q_2
\]
We use also  \cite[Theorem 8]{L.T} to deduce that
$(\nabla{u_k})$ is bounded in
${\textbf{L}}^{(\frac{2n}{n+2})^*,1}(B_1^n)={\textbf{L}}^{2,1}(B_1^n)$.
In particular,  there exist a constant $C > 0 $ depending only on
$n$ such that
\begin{equation}
\| \nabla {u_k}\|_{\textbf{L}^{2,1} (B_1^n)}\leq C \label{eq4.20}
\end{equation}
We deduce from (\ref{eq4.19}), (\ref{eq4.20}) together with the
${\textbf{L}}^{2,1}-{\textbf{L}}^{2,\infty}$ duality that
\[
\|\nabla{u_k}\|_{{\textbf{L}}^{2}(B_1^n\setminus{B_{R\delta_k}})}
\leq \|\nabla{u_k}\|_{{\textbf{L}}^{2,1}(B_1^n\setminus{B_{R\delta_k}})}
\|\nabla{u_k}\|_{{\textbf{L}}^{2,\infty}(B_1^n\setminus{B_{R\delta_k}})}\\
\leq C(\epsilon)
\]
for a constant $C(\varepsilon)\to 0$ as $\varepsilon\to 0$. Now,
we use the embedding
$\textbf{H}^1\hookrightarrow{\textbf{L}}^{2n/(n-2)}$
continuously, we obtain
\begin{align*}
\|u_k\|_{{\textbf{L}}^{2n/(n-2)}(B_1^n\setminus{B_{R\delta_k}})}
&\leq  C
\|\nabla{u_k}\|_{{\textbf{L}}^{2}(B_1^n\setminus{B_{R\delta_k}})}\\
&\leq C(\varepsilon)\to 0 \quad\text{as } \varepsilon\to  0 .
\end{align*}
We deduce that
\[
\lim_{R\to\infty} \lim_{k\to\infty}
\int_{B_{1}^n(0)\setminus{B_{R\delta_{k}(y_{k})}}}(|{u}_{k}|^{2n/(n-2)}
+|\nabla{{u}_{k}}|^{2})(x)\,dx = 0
\]
This proves  Theorem \ref{thm1.3} in the case of one bubble.

The case of more than one bubble can be handled in a
very similar way and we just give few details for $m = 2$. The proof
starts the same until (\ref{eq5.24}) which cannot hold any more
otherwise we would have had one bubble only as it is
(\ref{eq5.24}) holds. It remains to show that: for any
$\varepsilon \ge 0 $, there are sufficiently large R $> 0 $ and a
sequence $r_i \to 0 $ such that for any $ R {\delta}_i \le
r_i \le 1/2$,
\begin{equation}
\begin{gathered}
 \lim_{R \to \infty} \lim_{i
\to \infty} {\int_{\{0\} \times{ B^n_{r_i}
\setminus {B^n_{R {\delta}_i}}(0)}}}({  {
\frac{1}{2} \vert{\nabla{v_i}} \vert}^{2} +\frac{n-2}{2n}  {\vert
{v_i} \vert }^{2n/(n-2)}  })\,dx = 0\,,
\\
 \lim_{i \to \infty}
{\int_{\{0\}\times{ B^n_{1/2}\setminus{B^n_{r_i}(0)}}}
(  \frac{1}{2} {\vert {\nabla{v_i}} \vert}^{2} +\frac{n-2}{2n}
{\vert {v_i} \vert } ^{ 2n/(n-2)} )\,dx= 0}
\end{gathered} \label{eq5.31}
\end{equation}
where $v_i$ is defined  by $v_{i}(y)  = {r_i}^{(n-2)/2} $
$u_{i} (r_i y)$ , $y \in {\mathbb{R}}^n $.

The proof of (\ref{eq5.31})  can be done exactly as the proof of
(\ref{eq5.24}), the case of 2 bubbles is then proved. To prove the
general case,  for any number $m\ge 2$, one can follow
exactly the same strategy.

\end{document}